\documentclass{article}             
\usepackage{amssymb}
\def\mR{{\mathbb R}}
\def\fg{{\mathfrak g}}
\def\fh{{\mathfrak h}}
\def\fm{{\mathfrak m}}

\newtheorem       {teor}{Theorem}

\newtheorem{prop} [teor]{Proposition}
\newtheorem{lemma}[teor]{Lemma}
\newtheorem{cor}  [teor]{Corollary}

\begin{document}   
\begin{center}
\Large{The existence of light-like homogeneous geodesics
in~homogeneous Lorentzian manifolds}
\bigskip

\large{Zden\v ek Du\v sek}
\end {center}
\bigskip

\begin{abstract} 
In previous papers, a fundamental affine method for studying
homogeneous geodesics was developed. Using this method and elementary
differential topology it was proved that any homogeneous affine manifold
and in particular any homogeneous pseudo-Riemannian manifold admits
a homogeneous geodesic through arbitrary point.
In the present paper this affine method is refined and adapted
to the pseudo-Riemannian case. Using this method and elementary topology
it is proved that any homogeneous Lorentzian manifold of even dimension
admits a light-like homogeneous geodesic.
The method is illustrated in detail with an example of the Lie group
of dimension 3 with an invariant metric,
which does not admit any light-like homogeneous geodesic.
\end{abstract}
\bigskip

\noindent
{\bf MSClassification:} {53B05, 53C22, 53C30, 53C50}

\section{Introduction}
Let $M$ be a pseudo-Riemannian manifold.
If there is a connected Lie group $G\subset I_0(M)$
which acts transitively on $M$ as a group of isometries,
then $M$ is called a {\it homogeneous pseudo-Riemannian manifold\/}.
It can be naturally identified with the
{\it pseudo-Riemannian homogeneous space\/} $(G/H,g)$, where $H$ is the isotropy
group of the origin $p\in M$.

If the metric $g$ is positive definite, then
$(G/H,g)$ is always a {\it reductive homogeneous space\/}:
We denote by $\fg$ and $\fh$ the Lie algebras of $G$ and $H$ respectively
and consider the adjoint representation ${\rm Ad}\colon H\times\fg\rightarrow\fg$
of $H$ on $\fg$.
There exists the {\it reductive decomposition} of the form
$\fg=\fm+\fh$ where $\fm\subset\fg$
is a vector subspace such that ${\rm Ad}(H)(\fm)\subset\fm$.
For a fixed reductive decomposition $\fg=\fm+\fh$ there is the natural identification
of $\fm\subset\fg=T_eG$ with the tangent space $T_pM$ via the projection
$\pi\colon G\rightarrow G/H=M$.
Using this natural identification and the scalar product $g_p$ on $T_pM$,
we obtain the invariant scalar product $\langle\, ,\rangle $ on $\fm$.

If the metric $g$ is indefinite, the reductive decomposition may not exist
(see for instance \cite{FMP} or \cite{FR} for examples of nonreductive
pseudo-Riemannian homogeneous spaces).
In such a case, we can study the manifold $M$ using a more fundamental
affine method, which was proposed in \cite{DKV} and \cite{D3}.
It is based on the well known fact that homogeneous pseudo-Riemannian
manifold $M$ with the origin $p$ admits $n={\mathrm{dim}}M$
Killing vector fields which are linearly
independent at each point of some neighbourhood of $p$.

A geodesic $\gamma(s)$ through the point $p$ is {\it homogeneous}
if it is an orbit of a one-parameter group of isometries.
More explicitly, if $s$ is an affine parameter and $\gamma(s)$ is
defined in an open interval $J$, there exists
a diffeomorphism $s=\varphi(t)$ between the real line and
the open interval $J$ and
a nonzero vector $X\in\fg$ such that
$\gamma(\varphi(t))={\rm exp}(tX)(p)$ for all $t\in\mR$.
The vector $X$ is called {\it geodesic vector}.
The diffeomorphism $\varphi(t)$ may be nontrivial only for null curves
in a properly pseudo-Riemannian manifold.

In the reductive case, geodesic vectors are characterized
by the following {\it geodesic lemma}
(see \cite{KVa} for the Riemannian version, \cite{FMP} for the first formulation
in the pseudo-Riemannian case and \cite{DKa} for the complete mathematical proof).
\begin{lemma} \label{golema}
Let $X\in{\fg}$.
Then the curve $\gamma(t)={\rm exp}(tX)(p)$
is geodesic with respect to some parameter $s$ if and only if
\begin{eqnarray}
\nonumber
\langle [X,Z]_{\mathfrak m},X_{\mathfrak m}\rangle & = &
k\langle X_{\fm},Z\rangle
\end{eqnarray}
for all $Z\in{\mathfrak m}$ and for some constant $k\in{\mathbb{R}}$.
If $k=0$, then $t$ is an affine parameter for this geodesic.
If $k\neq 0$, then $s=e^{-kt}$ is an affine parameter for the geodesic.
The second case can occur only if the curve $\gamma(t)$ is a null curve
in a properly pseudo-Riemannian space.
\end{lemma}
In the paper \cite{KS}, it was proved that any homogeneous Riemannian manifold
admits a homogeneous geodesic through the origin.
The generalization to the pseudo-Riemannian (reductive and nonreductive) case
was obtained in \cite{D} in the framework of a more general result,
which says that any homogeneous affine manifold $(M,\nabla)$ admits a homogeneous
geodesic through the origin.
Here the affine method from \cite{DKV} and \cite{D3},
based on the study of integral curves of Killing vector fields, was used.
The proof is using differential topology, namely the degree
of a smooth mapping $S^{n}\rightarrow S^n$ without fixed points.

A homogeneous pseudo-Riemannian manifold all of whose geodesics are homogeneous
is called a pseudo-Riemannian {\it g.o.\,ma\-ni\-fold\/} or {\it g.o. space}.
Their analogues with noncompact isotropy group are {\it almost g.o. spaces}.
For many results and further references on homogeneous geodesics
in the reductive case see for example the survey paper \cite{D0}.

In pseudo-Riemannian geometry, null homogeneous geodesics are of particular interest.
In \cite{FMP} and \cite{P}, plane-wave limits (Penrose limits)
of homogeneous spacetimes along light-like homogeneous geodesics were studied.
However, it was not known whether any homogeneous pseudo-Riemannian or Lorentzian
manifold admits a null homogeneous geodesic.

In \cite{CM}, an example of a 3-dimensional Lie group with an invariant Lorentzian
metric which does not admit light-like homogeneous geodesic was described.
Here the standard geodesic lemma was used, because the example is reductive.

In the present paper, the affine method used in \cite{D}, \cite{D3} and \cite{DKV}
for the study of homogeneous affine manifolds is adapted
to the pseudo-Riemannian case.
As the main result it is shown that any Lorentzian homogeneous manifold
of even dimension admits a light-like homogeneous geodesic through the origin.
The~calculation is particularly easy in the case of a Lie group $G=M$
with a~left-invariant metric. As an illustration,
the method is applied on an example of a~Lie group from \cite{CM}.

\section{The main result}
Let $(M,g)$ be a homogeneous pseudo-Riemannian manifold of dimension $n$,
let $G$ be a group of isometries acting transitively on $M$ and let $p\in M$.
Let $\nabla$ be the induced pseudo-Riemannian connection on $M$.
It is well known that there exist $n$ Killing vector fields $K_1,\dots,K_n$ on $M$
which are linearly independent at each point of some neighbourhood $U$ of $p$.
Let $B=\{K_1(p),\dots,K_n(p)\}$ be the basis of the tangent space $T_pM$.
Any tangent vector $X\in T_pM$ has coordinates $(x^1,\dots,x^n)$ with respect
to the basis $B$ and it determines the Killing vector field
$X^*=x^1K_1+\dots+x^nK_n$ and the integral curve $\gamma_X$ of $X^*$ through $p$.
The following Proposition is a standard one.
\begin{prop}
Let $\phi_X(t)$ be the $1$-parameter group of isometries corresponding
to the Killing vector field $X^*$.
For all $t\in\mR$, it holds
\[
\phi_X(t)(p)=\gamma_X(t), \qquad
\phi_X(t)_*(X^*_p)=X^*_{\gamma_X(t)}.
\]
\end{prop}
It is well known that the covariant derivative $\nabla_{X^*}X^*$
depends only on the values of the vector field $X^*$ along the curve $\gamma_X(t)$.
From the invariance of the metric $g$ and the connection $\nabla$
with respect to the group $G$, we obtain the following:
\begin{prop}
\label{p2}
Along the curve ${\gamma_X(t)}$, it holds for all $t\in\mR$
\begin{eqnarray}
\nonumber
g_p(X^*,X^*) & = & g_{\gamma_X(t)}(X^*_{\gamma_X(t)}, X^*_{\gamma_X(t)} ), \cr
\phi_X(t)_*(\nabla_{X^*}X^*\big |_p) & = & \nabla_{X^*}X^*\big |_{\gamma_X(t)}.
\end{eqnarray}
\end{prop}
Now we formulate the crucial feature.
\begin{prop}
\label{p3}
Let $(M,g)$ be a homogeneous Lorentzian manifold, $p\in M$ and
$X\in T_pM$. Then, along the curve $\gamma_X(t)$, it holds
\begin{eqnarray}
\nonumber
\nabla_{X^*} X^*\big |_{\gamma_X(t)}\in (X^*_{\gamma_X(t)})^\perp.
\end{eqnarray}
\end{prop}
{\it Proof.}
We use the basic property $\nabla g=0$ in the form
\begin{eqnarray}
\label{eq3}
\nabla_{X^*} g(X^*,X^*) = 2g(\nabla_{X^*} X^*,X^*).
\end{eqnarray}
According to Proposition \ref{p2}, the function $g(X^*,X^*)$ is constant
along $\gamma_X(t)$.
Hence, the left-hand side of the equality (\ref{eq3}) is zero
and the right-hand side gives the statement.
$~\hfill\square$

\begin{teor}
\label{t1}
Let $(M,g)$ be a homogeneous Lorentzian manifold of even dimension $n$
and let $p\in M$.
There exist a light-like vector
$X\in T_pM$ such that along the integral curve $\gamma_X(t)$
of the Killing vector field $X^*$ it holds
\begin{eqnarray}
\label{f1}
\nabla_{X^*} X^*\big |_{\gamma_X(t)} = k\cdot X^*_{\gamma_X(t)},
\end{eqnarray}
where $k\in\mR$ is some constant.
\end{teor}
{\it Proof.}
Let us choose the Killing vector fields $K_1,\dots K_n$
such that the vectors $K_1(p),\dots,K_n(p)$
form a pseudo-orthonormal basis of $T_pM$ with $K_n(p)$ timelike.
Again, any airthmetic vector $x=(x^1,\dots,x^n)\in \mR^n$
determines the Killing vector field
$ X^*=\sum_{i=1}^n x^iK_i$.
Using the identification of $x$ with $X^*_p$ we identify $\mR^n$ with $T_pM$.
There is the natural scalar product on $\mR^n$
which comes from the scalar product $g_p$ on $T_pM$ and this identification.
Let us consider arithmetic vectors of the form
$x=(\tilde x,1)$, where $\tilde x\in S^{n-2}\subset\mR^{n-1}$.
For the corresponding vector field $X^*$, we have
$g_p(X^*_p,X^*_p)=0$ and the vectors $\tilde x\in S^{n-2}$
determine light-like directions in $\mR^n\simeq T_pM$.

For each light-like vector $x=(\tilde x,1)\in\mR^n\simeq T_pM$, we denote
$Y_x=\nabla_{X^*}{X^*}\big |_p$.
With respect to the pseudo-orthonormal basis $B=\{K_1(p),\dots,K_n(p)\}$,
we denote the components of the vector $Y_x$ as $y(x)=(y^1,\dots,y^n)$.
Using Proposition \ref{p3}, we see that $y(x)\perp x$.
We define the new vector $t_x$ as
\begin{eqnarray}
\nonumber
t_x & = & y(x) - y^n \cdot x.
\end{eqnarray}
Because $x$ is light-like vector, it holds also $t_x\perp x$.
For the components of $t_x$, we have $t_x=(\tilde t_x,0)$,
where $\tilde t_x\in\mR^{n-1}$.
We easily see that $\tilde t_x\perp \tilde x$, with respect to the positive
scalar product on $\mR^{n-1}$ which is the restriction of the indefinite
scalar product on $\mR^n$.
The assignment $\tilde x\mapsto\tilde t_x$ defines a smooth tangent vector
field on the sphere $S^{n-2}$.
If $n$ is even, then according to the well known topological theorems,
this vector field must have a zero value.
In other words, there exist a vector $\tilde x\in S^{n-2}$ such that
for the corresponding vector $x=(\tilde x,1)$ it holds $t_x=0$.
For this vector $x$, it holds $y(x)=k\cdot x$ and formula (\ref{f1})
for the corresponding Killing vector field $X^*$ is satisfied at $t=0$.
Using Proposition \ref{p2}, we obtain the formula for all $t\in\mR$.
$~\hfill\square$

\begin{cor}
Let $(M,g)$ be a homogeneous Lorentzian manifold of even dimension $n$
and let $p\in M$. There exist a light-like homogeneous geodesic through~$p$.
\end{cor}
{\it Proof.}
We consider the vector $X\in T_pM$ which satisfies Theorem \ref{t1}.
The integral curve $\gamma_X(t)$ through $p$
of the corresponding Killing vector field $X^*$ is homogeneous geodesic.
$~\hfill\square$

\section{Invariant metric on a Lie group}
Let $M=G$ be a Lie group with a left-invariant pseudo-Rieamannian metric $g$
acting on itself by left translations and let $p=e$ be the identity.
For any tangent vector $X\in T_eM$ and the corresponding Killing vector field $X^*$,
we consider the vector function $X^*_{\gamma_X(t)}$ along the integral curve
$\gamma_X(t)$ through $e$.
It can be uniquely extended to the left-invariant vector field $L^X$ on $G$.
Hence, along $\gamma_X$, we have
\begin{equation}
\label{lk}
L^X_{\gamma_X(t)}=X^*_{\gamma_X(t)}.
\end{equation}
At general points $q\in G$, values of left-invariant vector field $L_X$
do {\it not} coincide with the values of the Killing vector field $X^*$,
which is {\it right-invariant}.
However, as we are interested in calculations along the curve $\gamma_X(t)$,
we can work with respect to the moving frame of left-invariant vector fields
and use formula (\ref{lk}).
\begin{prop}
\label{p7}
Let $\{L_1,\dots,L_n\}$ be a left-invariant moving frame on a Lie group $G$
with a left-invariant pseudo-Riemannian metric $g$ and the induced
pseudo-Riemannian connection $\nabla$. Then it holds
\begin{eqnarray}
\nonumber
\nabla_{L_i}L_j  = \sum_{k=1}^n \gamma_{ij}^k L_k,\qquad i,j=1,\dots,n,
\end{eqnarray}
where $\gamma_{ij}^k$ are constants.
\end{prop}
{\it Proof.}
It follows from the invariance of the affine connection $\nabla$.
$~\hfill\square$

\medskip
Now we illustrate the affine method of the previous section with an example
of the 3-dimensional Lie group $E(1,1)$ with an invariant Lorentzian metric
which has no light-like homogeneous geodesic.
We choose one of the examples described in the paper \cite{CM}
by the standard method for reductive pseudo-Riemannian homogeneous manifolds
and the geodesic lemma.
We construct explicitly the vector field $\tilde t_x$, which has no zero value
in this case.

The group $E(1,1)$ can be represented by the matrices of the form
\begin{eqnarray}
\nonumber
\left (
\begin{array}{ccc}
e^{-w} & 0   & u \cr
0      & e^w & v \cr
0      & 0   & 1
\end{array}
\right ).
\end{eqnarray}
Hence, the manifold $M=E(1,1)$ can be identified with the 3-space $\mR^3[u,v,w]$.
The left-invariant vector fields are
$U=e^{-w}\partial_u,\, V=e^{w}\partial_v,\, W=\partial_w$.
We choose the new moving frame $\{E_1,E_2,E_3\}$ given as
\begin{eqnarray}
\nonumber
E_1  =  U - V,\qquad
E_2  =  - W,\qquad
E_3  =  1/2(U + V).
\end{eqnarray}
In this frame, we have the following rules for the Lie bracket
\begin{eqnarray}
\nonumber
[E_1,E_3]  =  0,\qquad
[E_2,E_1]  =  2 E_3,\qquad
[E_2,E_3]  =  1/2 E_1.
\end{eqnarray}
We introduce the pseudo-Riemannian metric $g$ such that the basis
determined by the above frame at any point $p\in M$
is pseudo-orthonormal basis of $T_pM$ with $E_3$ timelike
(we keep the notation from \cite{CM} here).

It is straightforward to write down the above metric $g$ in coordinates
in the form
\begin{eqnarray}
\nonumber
ds^2 & = & -\frac{1}{4} ( 3e^{2w} du^2
                        + 3e^{-2w} dv^2 + 10 du dv - 4 dw^2 )
\end{eqnarray}
and to calculate the nonzero Christoffel symbols
\begin{eqnarray}
\nonumber
\Gamma_{11}^{3} = \frac{3}{4} e^{2w}, &
\Gamma_{13}^{1} =-\frac{9}{16}, &
\Gamma_{13}^{2} = \frac{15}{16} e^{2w}, \cr
\Gamma_{22}^{3} =-\frac{3}{4} e^{-2w}, &
\Gamma_{23}^{2} = \frac{9}{16}, &
\Gamma_{23}^{1} =-\frac{15}{16} e^{-2w}.
\end{eqnarray}
However, we can write down the same information in the frame $\{E_1,E_2,E_3\}$.
By definition, we have at any point $p\in M$
\begin{eqnarray}
\nonumber
g(E_1,E_1)  =   g(E_2,E_2) = 1, \quad
g(E_3,E_3)  =   -1, \quad
g(E_i,E_j)  =   0, \quad i\neq j.
\end{eqnarray}
By the straightforward calculations, we obtain
nonzero covariant derivatives which satisfy Proposition \ref{p7}:
\begin{eqnarray}
\label{gijk}
\nabla_{E_1}E_2  = -\frac{3}{4}E_3 , \quad
\nabla_{E_1}E_3  = -\frac{3}{4}E_2 , \quad
\nabla_{E_2}E_3  =  \frac{5}{4}E_1 , \cr
\nabla_{E_2}E_1  =  \frac{5}{4}E_3 , \quad
\nabla_{E_3}E_1  = -\frac{3}{4}E_2 , \quad
\nabla_{E_3}E_2  =  \frac{3}{4}E_1 .
\end{eqnarray}
We will perform all calculations in this moving frame, or with respect
to the corresponding pseudo-orthonormal basis $B=\{ E_1(e),E_2(e),E_3(e)\}$
of the tangent space $T_eM\simeq\mR^3$ at the origin $e\in E(1,1)$.
Any arithmetic vector $x=(x^1,x^2,x^3)\in \mR^3$ determines
the left-invariant vector field
\begin{eqnarray}
\nonumber
L^X=x^1E_1+x^2E_2+x^3E_3.
\end{eqnarray}
We are interested in light-like vectors $X\in T_eM$,
hence $x=(\sin(\varphi),\cos(\varphi),1)$,
$\tilde x=(\sin(\varphi),\cos(\varphi) )\in S^1$ for some $\varphi\in\mR$.
For the corresponding left-invariant vector field $L^X$
we calculate using (\ref{gijk}) the covariant derivative
\begin{eqnarray}
\nonumber
\nabla_{L^X}{L^X} & = &
  2 \cos(\varphi) E_1 
-\frac{3}{2} \sin(\varphi) E_2 
+\frac{1}{2} \sin(\varphi)\cos(\varphi) E_3, \cr
y(x) & =  &  \Bigl ( 2 \cos(\varphi) ,
-\frac{3}{2} \sin(\varphi) ,
 \frac{1}{2} \sin(\varphi)\cos(\varphi) \Bigr ).
\end{eqnarray}
We see immediately that $y(x)\perp x$.
The projection $t_x$ is
\begin{eqnarray}
\nonumber
t_x& = & y(x) - \frac{1}{2} \sin(\varphi)\cos(\varphi)\cdot x= \cr
 & = & \Bigl ( 2\cos(\varphi)-\frac{1}{2} \sin^2(\varphi) \cos(\varphi) ,
               -\frac{3}{2}\sin(\varphi)-\frac{1}{2}\sin(\varphi)\cos^2(\varphi) ,
               0 \Bigr ) = \cr
 & = & \Bigl [ 2-\frac{1}{2}\sin^2(\varphi) \Bigr ] \cdot
       \Bigl ( \cos(\varphi) ,-\sin(\varphi), 0 \Bigr ).
\end{eqnarray}
We see that $t_x\perp x$ and $\tilde t_x\perp\tilde x$.
Clearly, $\tilde x\mapsto\tilde t_x$ defines the smooth vector field on $S^1$,
which is nonzero everywhere.
Hence, there is not any vector $X\in T_eG$ which satisfies Theorem \ref{t1}.

\begin{center}
{\bf Acknowledgements}
\end{center}
The author was supported by the grant GA\v{C}R 201/11/0356.

\bigskip
 
\noindent
Address of the author:\\
Zden\v ek Du\v sek\\
Palacky University, Faculty of Science\\
17. listopadu 12, 771 46 Olomouc, Czech Republic\\
zdenek.dusek@upol.cz\\
 
\end{document}